\documentclass[12pt]{article}
\usepackage{amsxtra,amssymb,amsthm,amsmath,latexsym}

\theoremstyle{plain}
\newtheorem{theorem}{Theorem}[section]

\newtheorem{remark}[theorem]{Remark}

\numberwithin{equation}{section}

\newcommand{\refT}[1]{Theorem~\ref{T:#1}}
\newcommand{\refS}[1]{Section~\ref{S:#1}}

\newcommand{\refR}[1]{Remark~\ref{R:#1}}

\def\lra{\longrightarrow}
\def\R{{\mathbb R}}

\def\ve{{\varepsilon}}

\def\circCinfty{{\overset{\circ}{C}\kern-.01in{}^\infty}}

\def\oH1{{\overset{\circ}{H}\kern-.04in{}^1}}

\def\bee{\begin{equation*}}
\def\eee{\end{equation*}}
\def\be{\begin{equation}}
\def\ee{\end{equation}}

\begin{document}
\title{A singular perturbation problem}

\author{
A.G. Ramm\\
 Mathematics Department, Kansas State University, \\
 Manhattan, KS 66506-2602, USA\\
ramm@math.ksu.edu\\
}
\date{}
\maketitle\thispagestyle{empty}

\begin{abstract}
\footnote{Math subject classification: 35J60, 35B25 }
\footnote{key words:  nonlinear elliptic equations, singular perturbations }

Consider the equation $-\ve^2\Delta u_\ve+q(x)u_\ve=f(u_\ve)$ in $\R^3$,
$|u(\infty)|<\infty$, $\ve=const>0$. Under what assumptions on $q(x)$ and
$f(u)$ can one prove that the solution $u_\ve$ exists and
$\lim_{\ve\to 0} u_\ve=u(x)$, where $u(x)$ solves the limiting problem
$q(x)u=f(u)$?
These are the questions discussed in the paper.
\end{abstract}


\section{Introduction}\label{S:1}
Let
\be\label{e1.1}
-\ve^2\Delta u_\ve+q(x)u_\ve=f(u_\ve) \hbox{\ in\ }\R^3,
\qquad |u_\ve(\infty)|<\infty, \ee
$\ve=const>0$, $f$ is a nonlinear smooth function, $q(x)\in C(\R^3)$ is a
real-valued function
\be\label{e1.2}
a^2\leq q(x),\qquad a=const>0. \ee
We are interested in the following questions: 

1)Under what assumptions does
problem \eqref{e1.1} have a solution? 

2)When does $u_\ve$ converge to $u$
as $\ve\to 0$? 

Here $u$ is a solution to
\be\label{e1.3}
q(x)u=f(u). \ee

The following is an answer to the first question.

\begin{theorem}\label{T:1.1}
Assume $q\in C(\R^3)$, (1.2) holds, 
$f(0)\not= 0$, and $a$ is sufficiently large (see \eqref{e2.7}
and \eqref{e2.9} below). Then equation \eqref{e1.1} has a solution 
$u_\ve \neq 0$,  $u_\ve \in C(\R^3),$ for any
$\ve>0$.
\end{theorem}

In Section 4 the potential $q$ is allowed to grow at infinity.

An answer to the second question is:

\begin{theorem}\label{T:1.2}
If $\frac{f(u)}{u} $ is a monotone, growing function, such 
that $\frac{f(u)}{u}\to \infty$ and
$\min_{u\geq u_0}\frac{f(u)}{u}< a^2$, 
where $u_0>0$ is a fixed number,
then there is a solution $u_\ve$ to \eqref{e1.1} such that
\be\label{e1.4}
\lim_{\ve\to 0} u_\ve(x)=u(x), \ee
where $u(x)$ solves \eqref{e1.3}.
\end{theorem}

Singular perturbation problems have been discussed in the 
literature  \cite{B}, \cite{L},  \cite{VL}.

In \refS{2} proofs are given.

In \refS{3} an alternative approach is proposed. 

In \refS{4} an extension of the results to a larger class of potentials
is given.

\section{Proofs}\label{S:2}

\begin{proof}[Proof of \refT{1.1}]
The existence of a solution to (1.1) is proved by means of the
contraction mapping principle.

Let $g$ be the Green's function
\be\label{e2.1}
(-\ve^2\Delta+a^2)g=\delta(x-y)\hbox{\ in\ }\R^3, 
\qquad g:=g_a(x,y,\ve)\underset{|x|\to\infty}{\lra} 0,
\qquad g=\frac{e^{-\frac{a}{\ve}|x-y|}}{4\pi|x-y|\ve^2}. \ee
Let $p:=q-a^2\geq 0$. Then \eqref{e1.1} can be written as:
\be\label{e2.2}
u_\ve(x)=-\int_{\R^3} gpu_\ve dy+\int_{\R^3} gf(u_\ve)dy:=T(u_\ve). \ee
Let $X=C(\R^3)$ be the Banach space of continuous and globally bounded 
functions, 
 $B_R:=\{ v:\|v\|\leq R\}$, and $\|v\|:=\sup_{x\in\R^3} |v(x)|$.

We choose $R$ such that
\be\label{e2.3}
T(B_R)\subset B_R \ee
and
\be\label{e2.4}
\|T(v)-T(w)\|\leq \gamma\|v-w\|, \qquad v,w\in B_R, \qquad 0<\gamma<1.\ee
If \eqref{e2.3} and \eqref{e2.4} hold, then the contraction mapping principle 
yields a solution $u_\ve\in B_R$ to \eqref{e2.2}, and, therefore, to 
problem (1.1).

The assumption $f(0)\not= 0$ guarantees that $u_\ve\not=0$.

Let us check \eqref{e2.3}. If $\|v\|\leq R$, then 
\be\label{e2.5}
\|T(v)\|\leq \|v\| \|p\| \|\int_{\R^3} g(x,y)dy\|
+\frac{M(R)}{a^2} \leq \frac{\|p\|R+M(R)}{a^2} \leq R,\ee
where $M(R):=\max_{|u|\leq R} |f(u)|$. 
Here we have used the following estimate:
\be\label{e2.6}
\int_{\R^3} g(x,y)dy=\int_{\R^3} 
\frac{e^{-\frac{a}{\ve}|x-y|}}{4\pi|x-y|\ve^2} \ dy=\frac{1}{a^2}. \ee
If $\|p\|<\infty$ and $a$ is such that
\be\label{e2.7}
\frac{\|p\|R + M(R)}{a^2} \leq R,\ee
then \eqref{e2.3} holds.

Let us check \eqref{e2.4}. Assume that $v,w\in B_R$, $v-w:=z$. Then 
\be\label{e2.8}
\|T(v)-T(w)\| \leq \frac{\|p\|}{a^2} \|z\|+ \frac{M_1(R)}{a^2} \|z\|, \ee
where $M_1(R)=\max_{\substack{|u|\leq R\\ |w|\leq R\\0\leq s\leq 1}}
|f'(u+sw)| \leq \max_{|\xi|\leq 2R} |f'(\xi)|$.
If 
\be\label{e2.9}
\frac{\|p\|+M_1(R)}{a^2} \leq \gamma< 1, \ee
then \eqref{e2.4} holds. By the contraction mapping principle, 
\eqref{e2.7}
and \eqref{e2.9} imply the existence and uniqueness of the solution
$u_\ve(x)$ to \eqref{e1.1} in $B_R$ for any $\ve>0$.

\refT{1.1} is proved.
\end{proof}

\begin{proof}[Proof of \refT{1.2}]
In the proof of \refT{1.1} one can choose $R$ and $\gamma$ independent
of $\ve>0$. Let us denote by $T_{\ve}$ the operator defined in (2.2).
Then (see Remark 2.2 ) one 
has
\be\label{e2.10}
\lim_{\ve\to 0}\|T_\ve(v)-T_0(v)\|=0, \ee
where
\be\label{e2.11}
T_0(v)=\frac{-pv+f(v)}{a^2}.\ee
It is known \cite{KA} and easy to prove (see \refR{2.3}) that if
\eqref{e2.10} holds for every $v\in X$, and $\gamma$ in \eqref{e2.4}
does not depend on $\ve$, then \eqref{e1.4} holds, where $u$ solves
the limiting equation \eqref{e2.2}:
\be\label{e2.12}
u=T_0(u)=\frac{-pu+f(u)}{a^2}. \ee
Equation \eqref{e2.12} is equivalent to \eqref{e1.3}.
Theorem 1.2 is proved.
\end{proof}

\begin{remark}\label{R:2.1}
Conditions of \refT{1.1} and \refT{1.2} are satisfied if, for example,
$q(x)=a^2+1+\sin(\omega x)$,where $\omega =const>0$, 
$f(u)=(u+1)^m$, $m>1,$ 
or $f(u)=e^u$.
\end{remark}

\begin{remark}\label{R:2.2}
Note that in the distribution sense
\be\label{e2.13}
g_a(x,y,\ve)\to\frac{1}{a^2} \delta(x-y), \qquad \ve\to 0.\ee
\end{remark}

\begin{remark}\label{R:2.3}
Let $u=T_\ve(u)$, $v=T_{\ve_0}(v):=T_0(v)$,
and $T_\ve(w)\to T_0(w)$ for all $w\in X$,
$\|T_\ve(v)-T_\ve(w)\|\leq \gamma\|v-w\|$,
$0<\gamma<1$, $\gamma$ does not depend on $\ve$,
$u_{n+1}=T_\ve(u_n)$, $u_0=v$.
Then $u_1=T_\ve v$, and $\|u_n-v\|\leq\frac{1}{1-\gamma}\|u_1-v\|$.
 Taking $n\to\infty$, one gets
$\|u-v\|\leq\frac{1}{1-\gamma}\|T_\ve(v)-T_0(v)\|\to 0$ as $\ve\to\ve_0$.
\end{remark}

\section{A different approach}\label{S:3}
Let us outline a different approach to problem \eqref{e1.1}.
Set $x=\xi+\ve y$.
Then
\be\label{e3.1}
-\Delta_y w_\ve+a^2w_\ve+p(\ve y+\xi)w_\ve =f(w_\ve),
\qquad |w_\ve(\infty)|<\infty, \ee
$w_\ve:=u_\ve(\ve y+\xi)$, $p:=q(\ve y+\xi)-a^2\geq 0$. 
Thus
\be\label{e3.2}
w_\ve=-\int_{\R^3} G(x,y) p(\ve y+\xi) w_\ve dy
+\int_{\R^3} G(x,y) f(w_\ve)dy, \ee
where
\be\label{e3.3}
(-\Delta+a^2) G=\delta(x-y)\hbox{\ in\ }\R^3, 
\qquad G=\frac{e^{- a|x-y|}}{4\pi|x-y|}, \quad a>0.\ee
One has
\be\label{e3.4}
\int_{\R^3} G(x,y)dy=\frac{1}{a^2}. \ee
Using an argument similar to the one in the proofs of \refT{1.1} and
\refT{1.2}, one concludes that for any $\ve >0$  and any 
sufficiently large 
$a$, problem \eqref{e3.1} has a unique solution, which tends 
to a limit $w=w(y,\xi)$ as $\ve\to0$, where $w$ solves the problem
\be\label{e3.5}
-\Delta_y w+q(\xi)w=f(w), \qquad |w(\infty,\xi)|<\infty. \ee
Problem \eqref{e3.5} has an obvious solution $w=w(\xi)$, which is indepent
of $y$ and solves the equation
\be\label{e3.6}
q(\xi)w=f(w).\ee
The solution to \eqref{e3.5} is unique if $a$ is sufficiently large.
This is proved similarly to the proof of (2.9). Namely, let 
$b^2:=q(\xi)$. Note that $b\geq a$. If 
there are two solutions to (3.5), say $w$ and $v$, and if $z:=w-v$,
then $||z||\leq b^{-2}M_1(R)||z||<||z||$, provided that $b^{-2}M_1(R)<1$.
Thus $z=0$, and the uniqueness of the solution to \eqref{e3.5} is proved.

Replacing $\xi$ by $x$ in \eqref{e3.6}, we obtain the solution
found in \refT{1.2}.

\section{Extension of the results to a larger class of potentials}\label{S:4}
Here a method for a study of problem \eqref{e1.1} for a larger class of
potentials $q(x)$ is given. We assume that $q(x)\geq a^2$ and can grow
to infinity as $|x|\to\infty$. Note that in Sections 1 and 2 the potential
was assumed to be a bounded function. 
Let $g_\ve$ be the Green's function
\be\label{e4.1}
-\ve^2\Delta g_\ve + q(x)g_\ve=\delta(x-y) \hbox{\ in\ }\R^3,
\qquad |g_\ve(\infty,y)|<\infty.\ee
As in \refS{2}, problem \eqref{e1.1} is equivalent to 
\be\label{e4.2}
u_\ve=\int_{\R^3} g_\ve f(u_\ve(y))dy,\ee
and this equation has a unique solution in $B_R$ if $a^2$ is
sufficiently large. The proof, similar to the one given in \refS{2},
requires the estimate
\be\label{e4.3}
\int_{\R^3}g_\ve(x,y)dy\leq \frac{1}{a^2}.
\ee
Let us prove the above inequality.
Let $G_j$ be the Green's function satisfying equation (4.1) with $q=q_j$,
$j=1,2$.
Estimate (4.3) follows from the inequality 
\be\label{e4.4}
G_1\leq G_2 \quad \text { if } q_1\geq q_2.
\ee
This inequality can be derived from the maximum principle. 

If  $q_2=a^2$,
then $G_2=\frac{e^{-\frac{a}{\ve}|x-y|}}{4\pi|x-y|\ve^2}$, and the
inequality $g_\ve(x,y)\leq \frac{e^{-\frac{a}{\ve}|x-y|}}{4\pi|x-y|\ve^2}$
implies (4.3).

We prove below the following relation:
\be\label{e4.5}
\lim_{\ve\to 0} \int_{\R^3} g_\ve (x,y) h(y)dy=\frac{h(x)}{q(x)}
\qquad \forall h\in \circCinfty(\R^3), \ee
where $\circCinfty$ is the set of $C^\infty(\R^3)$ functions
vanishing at infinity together with their derivatives.
This formula is an analog to \eqref{e2.13}.

To prove \eqref{e4.5}, multiply \eqref{e4.1} by $h(y)$,
integrate over $\R^3$ with respect to $y$, then integrate
the first term by parts, and then let $\ve\to 0$. The result
is \eqref{e4.5}. 

Thus, \refT{1.1} and \refT{1.2} remain
valid for $q(x)\geq a^2$, $a>0$ sufficiently large,
$\frac{f(u)}{u}$ monotonically growing to infinity, and the
solution $u(x)$ to the limiting equation \eqref{e1.3} is the
limit of the solution to \eqref{e4.2} as $\ve\to 0$.

\end{document}